\documentclass[oneside,12pt]{article}
\usepackage[colorlinks,citecolor=DarkRed,urlcolor=DarkRed]{hyperref} 
\usepackage[x11names,svgnames]{xcolor} 
\usepackage{colorweb}                  
\usepackage{tocloft}                   
\usepackage{setspace}                  
\usepackage{indentfirst}               
\usepackage[a4paper
           ]{geometry}
\usepackage{edrx26a}                     
\def\dednatjobname{2026bad-foundations}  
\usepackage[backend=biber,
   style=alphabetic]{biblatex}    
\begin{document}

\ifluatex
  \catcode`\^^J=10
  \directlua{dofile "dednat7load.lua"}  
  \directlua{dednat7preamble()}         
  \directlua{dednat7oldheads()}         
\else
  \input\dednatjobname.dnt              
  \def\directlua#1{}
  \def\pu{}
\fi
\def\subfilepu{}

\def\saexpr#1{%
  \directlua{output:saexpr([==[#1]==])}%
  \ga{#1}%
  }

\def\Caurl   #1{\expr{Caurl("#1")}}    
\def\Cahref#1#2{\href{\Caurl{#1}}{#2}}
\def\Ca      #1{\Cahref{#1}{#1}}

\subfilepu

\def\pictgridstyle{\color{GrayPale}\linethickness{0.3pt}}
\def\pictaxesstyle{\linethickness{0.5pt}}
\def\pictnaxesstyle{\color{GrayPale}\linethickness{0.5pt}}
\celllower=2.5pt

\subfilepu
\def\picturedots(#1,#2)(#3,#4)#5{\saexpr{
    PictureDots.from(#1,#2, #3,#4, "#5"):topict():pgat("patc")
  }}

\def\picturedotsadef#1(#2,#3)(#4,#5)#6{\saexpr{
      PictureDots.from(#2,#3, #4,#5, "#6"):topict():pgat("patc")
  }}
\def\picturedotsdef #1(#2,#3)(#4,#5)#6{\saexpr{
      PictureDots.from(#2,#3, #4,#5, "#6"):topict():pgat("ptc")
  }}

\subfilepu

\subfilepu
\def\vbtc       #1{\myvcenter{\vbt{#1}}}
\def\vbtsmash   #1{\nodepthnoheight{\vbtc{#1}}}
\def\vbtbgboxcolor{YellowOrangeLight!45}

\def\co#1{{%
  \def\%{\char37}%
  \def\\{\char92}%
  \def\^{\char94}%
  \def\~{\char126}%
  \tt#1%
  }}
\def\qco#1{`\co{#1}'}
\def\qqco#1{``\co{#1}''}

\sa{tictactoe-.}{\phantom{\textsf{X}}}
\sa{tictactoe-x}{\textsf{X}}
\sa{tictactoe-o}{\textsf{O}}
\sa{tictactoe-X}{\ColorRed{\textsf{X}}}
\sa{tictactoe-O}{\ColorRed{\textsf{O}}}
\def\tictactoe#1#2#3#4#5#6#7#8#9{
  \ensuremath{
    \setlength{\arraycolsep}{2pt}
    \begin{array}[c]{c|c|c}
      \ga{tictactoe-#1} & \ga{tictactoe-#2} & \ga{tictactoe-#3} \\\hline
      \ga{tictactoe-#4} & \ga{tictactoe-#5} & \ga{tictactoe-#6} \\\hline
      \ga{tictactoe-#7} & \ga{tictactoe-#8} & \ga{tictactoe-#9} \\
    \end{array}
 }}

%
\subfilepu
\def\linktopage       #1#2{\hyperlink{page.#1}{#2}}
\def\mytocslide     #1#2#3{\par    \linktopage{#1}{#3} \dotfill\;\linktopage{#3}{#3}}
\def\toclineslidetex  #1#2{\par    \linktopage{#2}{#1} \dotfill \linktopage{#2}{#2}}
\def\toclineslidenttex#1#2{\par\;\;\linktopage{#2}{#1} \dotfill \linktopage{#2}{#2}}
\def\toclineslidelua  #1#2{\directlua{
    toclines:add {style="slide",   body=[=[#1]=], page=#2}
  }}
\def\toclineslidentlua#1#2{\directlua{
    toclines:add {style="slident", body=[=[#1]=], page=#2}
  }}
\def\SLIDE  #1{\par{\bf #1}\par \toclineslidelua  {#1}{\thepage}}
\def\SLIDENT#1{                 \toclineslidentlua{#1}{\thepage}}
\def\addslide#1#2{\SLIDE{#2}}

\def\Rq{\ColorRed{?}}



\pu

%
%
%

\def\comprehensionbox#1{\vbtbgbox{\ensuremath{\mat{#1}}}}
\def\V{\mathbf{V}}
\def\F{\mathbf{F}}
\def\Stop{\omit\vrule\phantom{$\scriptstyle($}\hss}
\def\Stop{\omit|\hss}
\def\HLine{\hline\\[0pt]}
\def\HLine{\hline}

\def\undtext#1#2{\underbrace{\mathstrut#1}_{\mathstrut\text{#2}}}
\def\ug#1{\undtext{#1}{gen}}
\def\uf#1{\undtext{#1}{filt}}
\def\ue#1{\undtext{#1}{expr}}


\sa {compr 10*a}      { \setofst     {10a}      {a∈\{2,3,4\}}             }
\sa {compr 10*a u}    { \setofst {\ue{10a}} {\ug{a∈\{2,3,4\}}}            }
\sa {compr 10*a u;}   { \setofsc            {\ug{a∈\{2,3,4\}}} {\ue{10a}} }

\sa {compr a^2<10}    { \setofst     {a∈\{2,3,4\}}      {a^2<10}           }
\sa {compr a^2<10 u}  { \setofst {\ug{a∈\{2,3,4\}}} {\uf{a^2<10}}          }
\sa {compr a^2<10 u;} { \setofsc {\ug{a∈\{2,3,4\}},  \uf{a^2<10}} {\ue{a}} }

\sa {compr triangle;} {
  \setofsc {    x∈\{1,\ldots,5\},
                y∈\{x,\ldots,6-x\}
           }
           {   (x,y) }
  }

\sa {compr triangle u;} {
  \setofsc {\ug{x∈\{1,\ldots,5\}},
            \ug{y∈\{x,\ldots,6-x\}}
           }
           {\ue{(x,y)}}
  }

\sa {compr triangle dots} {
  \picturedots(0,0)(6,6){ 1,1 1,2 1,3 1,4 1,5 2,2 2,3 2,4 3,3 }
  }

\sa {compr triangle dots wrong} {
  \picturedots(0,0)(6,6){ 1,1 1,2 1,3 1,4 1,5 2,2 2,3 2,4 3,3 4,3 }
  }

\sa {compr triangle tree thin} {
  \comprehensionbox{
    x & y & (x,y) \\\HLine
    1 & 1 & (1,1) \\
      & 2 & (1,2) \\
      & 3 & (1,3) \\
      & 4 & (1,4) \\
      & 5 & (1,5) \\
    2 & 2 & (2,2) \\
      & 3 & (2,3) \\
      & 4 & (2,4) \\
    3 & 3 & (3,3) \\
    4 & \Stop \\
    5 & \Stop \\
  }}

\sa {compr triangle tree wide} {
  \comprehensionbox{
    x & 6-x & \{x,\ldots,6-x\} & y & (x,y) \\\HLine
    1 & 5 & \{1,2,3,4,5\} & 1 & (1,1) \\
      &   &               & 2 & (1,2) \\
      &   &               & 3 & (1,3) \\
      &   &               & 4 & (1,4) \\
      &   &               & 5 & (1,5) \\
    2 & 4 & \{2,3,4\} & 2 & (2,2) \\
      &   &           & 3 & (2,3) \\
      &   &           & 4 & (2,4) \\
    3 & 3 & \{3\} & 3 & (3,3) \\
    4 & 2 & \{\} & \Stop \\
    5 & 1 & \{\} & \Stop \\
  }}


\def\comprehensionbox#1{\ensuremath{\sm{#1}}}
\def\comprehensionbox#1{\ensuremath{\mat{#1}}}
\def\comprehensionbox#1{\ensuremath{\fbox{\sm{#1}}}}
\def\comprehensionbox#1{\fbox{\ensuremath{\mat{#1}}}}
\def\comprehensionbox#1{\fbox{\ensuremath{\sm{#1}}}}
\def\comprehensionbox#1{\vbtbgbox{\ensuremath{\mat{#1}}}}
\def\V{\mathbf{V}}
\def\F{\mathbf{F}}
\def\Stop{\omit\vrule\phantom{$\scriptstyle($}\hss}
\def\Stop{\omit|\hss}
\def\HLine{\hline\\[0pt]}
\def\HLine{\hline}


\def\drafturl{http://anggtwu.net/LATEX/2025-2-C2.pdf}
\def\drafturl{http://anggtwu.net/2025.2-C2.html}
\def\drafturl{https://anggtwu.net/math-b.html\#2026-bad-foundations}
\def\draftfooter{\tiny \href{\drafturl}{\jobname{}} \ColorBrown{\shorttoday{} \hours}}

\def\rarr{\ColorRed{⇒}}


%
%


\title{Bad Foundations \\ and Manipulable Objects}

\author{%
  Eduardo Ochs%
  \thanks{eduardoochs@gmail.com \\
    \url{https://anggtwu.net/math-b.html\#2026-bad-foundations}} \\
  }

\maketitle


\begin{abstract}

  Imagine a student---let's call him `E', and make him a ``he''---that
  is enrolled in Calculus 2, and who believes that to pass in maths
  courses he only needs to memorize methods and apply them quickly and
  without errors. Let's imagine that `E' is an `E'xtreme case of ``bad
  foundations'' and that he knows how to solve $x+2=5$ by doing
  $x=5-2=3$, but he doesn't know how to substitute the $x$ in $x+2=5$
  by 3, and the only way that he knows of ``testing the solution'' is
  to apply the same method again and check that he got the same
  result.

  When we are teaching Calculus to classes that have many students
  that are extreme cases of bad foundations we need new strategies and
  tools; for example, we can't pretend that ``taking a particular
  case'' is an obvious operation anymore --- instead we need ways to
  make these operations easy to visualize. This article shows a way to
  do that using Maxima.

\end{abstract}


\def\Dedic{
  \begin{tabular}{l}
    Para Walter Machado Pinheiro, \\
    que não leu e não vai ler     \\
    documento nenhum, e se ler    \\
    não vai entender              \\
  \end{tabular}
}





\section{Introduction}
\label{introduction}

Let me start by some stories. I teach Calculus 2 and 3 in a
low-prestige campus of a high-prestige federal university in Brazil,
and all the characters that I am going to mention are based on real
students that I had in Calculus 2 -- that I will sometimes abbreviate
as ``C2''.


\msk

`B'en is a student that has `B'ad foundations and `E'ric is a student
who has `E'extremely bad foundations. They are taking Calculus 2 with
me, and in the second test, that has some questions abouts ODEs, they
commit the same error: they find $y=\sqrt{25-x^2}+3$ as the solution
of a certain ODE, and in the item in which they have to draw its
graph, they simplify $y=\sqrt{25-x^2}+3$ by doing this:
%
%
$$\begin{array}{rcl}
  y &\eqnp{1}& \sqrt{25-x^2} + 3 \\
    &=& \sqrt{5^2 - x^2} + 3 \\
    &=& 5 - x + 3 \\
    &=& 8 - x \\
  \end{array}
$$

I tell them that in the third `$=$' sign they used a rule that is
invalid. They panic, and they say: ``oh no, sorry, we got distracted
-- that won't happen again!!!''. Then I ask them to tell me what rule
they used there, and I discover that they don't have any idea of what
I mean... for them Maths is about memorizing methods and applying
them; ``understanding'' is only for geniuses
(\cite[p.7]{SchoenfeldWhenGood}), and all their knowledge about Maths
is {\sl procedural}, not {\sl conceptual}
(\cite{EngelbrechtBergsten}). I will explain more about bad
foundations on section \ref{bad-foundations}.

\msk


`Q'uentin is a student who believes that the goal of learning Maths is
to get answers `Q'uickly (\cite[p.33]{Boaler}). In a test question
that asked to calculate $\intx{\sqrt{1-x^2}}$ his answer was just
this:
$$\frac12 \left( \arcsin(x) + x\sqrt{1-x^2} \right)$$
I tried to explain to him that the ``right'' answer to that question
would a calculation with many steps, in which each `=' would be ``easy
to understand, to justify, and to verify'', but none of my arguments
made any sense to him.


\msk

`F'abian is a student who is more `F'riendly and talkative than most.
He tells me that he had very good grades in Calculus 1 and that he was
very good at calculating derivatives, and I decide that he will be my
reference for understanding what the students know about `F'unctions.
I ask him to show me how he would calculate
$\ddx \sin(\cos(\tan(42x))$; he tries but he gets stuck, and my hints
don't help him. It turns out that he learned the Chain Rule by a
method similar to the one in \cite[p.154]{Stewart8}, in which he had
to {\sl remember} what are the ``outer function'', the ``inner
function'', and their derivatives; he didn't have a way to write these
four functions down on paper.

The method the Fabian learned is optimized for speed, and in complex
cases, like $\ddx \sin(\cos(\tan(42x))$, it needs a lot of working
memory to be carried out... and Fabian had a memory overflow.

\msk


\sa{ddx sin(42x) diagram}{
  \sa{nw}{\pmat{\ddx f(g(x))    = \\ f'(g(x))g'(x)}}
  \sa{ne}{\pmat{\ddx f(42x)     = \\ f'(42x)·42}}
  \sa{sw}{\pmat{\ddx \sin(g(x)) = \\ \cos(g(x))g'(x)}}
  \sa{se}{\pmat{\ddx \sin(42x)  = \\ \cos(42x)·42}}
  \begin{array}{ccc}
    \ga{nw} & →        & \ga{ne} \\\\[-5pt]
      ↓     & \searrow &  ↓      \\\\[-5pt]
    \ga{sw} & →        & \ga{se} \\
  \end{array}
  }

{}

One day I ask all the students to calculate
$\ddx \,f(\sin(x^4) + \ln x)$. They all drop the `$f$', and the best
students get the result $\cos(x^4)·4x^3 + \frac{1}{x}$. I ask them how
they got there, and they say that the `$f$' was ``just notation'' and
meant ``function'', so they calculated ``the derivative of the
function $\sin(x^4) + \ln x$''.

I discuss that with Fabian, and I draw this diagram:
$$\ga{ddx sin(42x) diagram}$$

We name its nodes as $\sm{(1)&→&(2)\\↓&\searrow&↓\\(3)&→&(4)}$, and
Fabian tells me that (4) makes all sense to him and he recognizes (1)
as being the formula of the Chain Rule, but he doesn't have any idea
of what are (2) and (3).

%
This reveals that the students know far less than they should, but is
not totally surprising -- we know that books written in the 13th, 14th
and 15th centuries didn't treat unknowns like $x$ as entities with the
same algebraic properties as numbers (\cite{FilloyRojano},
\cite[p.216]{PuigRojano}), and the nodes (1), (2), and (3) in the
diagram above talk about unknown {\sl functions}, which are much worse
conceptually than unknown {\sl numbers}.

What did these students learn in Calculus 1? My colleagues share
practically nothing of their didactic strategies and teaching
materials, but they brag that in their courses they present
proofs, and that the students wouldn't see any proofs
if we let the engineers teach Calculus to them -- see
\cite{AhmadApplebyEdwards}...

Here's a hypothesis:

\begin{quote}

  When my students took C1 they saw a version of the Chain Rule ``for
  analysts'', in which (1) is a part of the statement of a theorem,
  and in the full statement $f$ and $g$ are quantified -- either as
  ``for any $f,g:\R→\R$ differentiable everywhere'', or as something
  like that but with different domains and different requirements of
  differentiability -- but the students retained very little from what
  they saw because they are in a ``procedural'' stage, and they didn't
  know how to operate on those quantifiers.

\end{quote}

Let me take a paragraph from \cite[p.62]{CarlsonGulick}:

\begin{quote}

  At the same time certain participants emphasized that the
  fundamental concepts in calculus must remain in the course, along
  with at least a certain amount of drill work to develop manipulative
  skill. We must minimize the tedium of working problems where no
  thought at all is necessary. A basic question is the following: What
  would be an appropriate blend between geometrically-motivated
  concepts, definitions and theorems, relevant applications, and
  rigorous proofs?

\end{quote}


\sa {subst early} {\begin{pmatrix}&{\frac{d}{d\,x}}\,f\left(g\left(x\right)\right)\cr \mbox{ = }&f'\left(g\left(x\right)\right)\,g'\left(x\right)\cr \end{pmatrix}\,\begin{bmatrix}g\left(x\right) := 42\,x\cr g'\left(x\right) := 42\cr \end{bmatrix}=\begin{pmatrix}&{\frac{d}{d\,x}}\,f\left(42\,x\right)\cr \mbox{ = }&f'\left(42\,x\right)\,42\cr \end{pmatrix}}

In this article I will discuss an unusual answer for the question
``What would be an appropriate blend...?'' above. My proposal is that
{\sl when we have lots of students like the ones that I described
  above} it is a good idea to present expressions and trees as some of
our most basic objects, and present this substitution operation
$$\ga{subst early}$$
very early in the course, before theorems and quantifiers.

\msk

\def\eqnp#1{\overset{\scriptscriptstyle{(#1)}}{=}}

\def\foo#1#2#3#4#5{ \text{#1} & #2 & =^#3 & #4 & \text{#5} \\ }
\def\foo#1#2#3#4#5{ \text{#1\;\;\;} & #2 & =^#3 & #4 & \text{\;\;\;#5} \\ }
\def\foo#1#2#3#4#5{ \text{#1\;\;\;} & #2 & \eqnp{#3} & #4 & \text{\;\;\;#5} \\ }
  
The figure below has a `L'eft part with particles `if', `then', and
`and', a `M'iddle part with equalities, and a `R'ight part with
justifications:
$$%
  \begin{array}{rrcll}
  \foo {If}          {f(g(x))} {1} {x}               {}
  \foo {then}   {\ddx f(g(x))} {2} {\ddx x}          {by (1)}
  \foo {and}                {} {3} {1}               {}
  \foo {and}    {\ddx f(g(x))} {4} {f'(g(x)) g'(x)}  {the chain rule}
  \foo {and}  {f'(g(x)) g'(x)} {5} {1}               {by (4), (2), (3)}
  \foo {and}           {g'(x)} {6} {1 / f'(g(x))}    {by (5)}
  \end{array}
$$

One day I write the part `M' on the blackboard and I ask the students
a) if that is how they learned that $\ddx \ln x = \frac1x$ in Calculus
1, and b) if all the steps are ``obvious'' to them. `ST'eve says that
that's not obvious at all, and asks me if I can make that clearer. I
show that we can write the columns `L' and `R', and I ask him what is
the first step that doesn't make sense to him. He answers that it's
`$\eqnp{5}$', and I write this expanded version:
$$%
  \begin{array}{rrcll}
  \foo {If}          {f(g(x))} {1} {x}               {}
  \foo {then}   {\ddx f(g(x))} {2} {\ddx x}          {by (1)}
  \foo {and}                {} {3} {1}               {}
  \foo {and}    {\ddx f(g(x))} {4} {f'(g(x)) g'(x)}  {the chain rule}
  \foo {and}  {f'(g(x)) g'(x)} {5} {\ddx f(g(x))}    {by (4)}
  \foo {and}                {} {6} {\ddx x}          {by (2)}
  \foo {and}                {} {7} {1}               {}
  \foo {and}  {f'(g(x)) g'(x)} {8} {1}               {by (5), (6), (7)}
  \foo {and}           {g'(x)} {9} {1 / f'(g(x))}    {by (5)}
  \end{array}
$$

%
We discover that `ST'eve did not know that the equality is `S'ymmetric
and `T'ransitive. He only knew three meanings for the `=' sign:
``compute'', ``simplify'', and ``solve this'' -- see for example
\cite[p.10]{FischbeinTacit} and \cite[p.179]{ThomasRethinking}.


\section{Expandable objects}
\label{expandable-objects}


{\bf Game trees.} The tree below is {\sl essentially the same} as the
one in \cite[Section 11.8]{Hutton},
$$\def\TTT{\tictactoe}
  \sa{tictactoe-O}{\ColorOrange{\textsf{O}}}
  \sa{tictactoe-X}{\ColorOrange{\textsf{X}}}
  \sa{tictactoe-X}{\ColorGreen{\textsf{X}}}
  \scalebox{0.75}{$\mat{
    \TTT o.. xxo xo. & \TTT oO. xxo xo. & \TTT ooX xXo Xo. \\\\[-6pt]
                     &                  & \TTT oo. xxo xoX & \TTT OOO xxo xox \\\\[-6pt]
                     & \TTT o.O xxo xo. & \TTT oXo xxo xo. & \TTT oxO xxO xoO \\\\[-6pt]
                     &                  & \TTT o.o xxo xoX & \TTT OOO xxo xoX \\\\[-6pt]
                     & \TTT o.. xxo xoO & \TTT oX. xxo xoo & \TTT oxO xxO xoO \\\\[-6pt]
                     &                  & \TTT o.X xXo Xoo \\
  }$}
$$
but the one in Hutton's book has the root node at the top, has lines
from each parent node to its child nodes, and is in black and white;
our tree draws child nodes to the right, like this,
$\bsm{1 & 1.1 & 1.1.1 \\ && 1.1.2 \\ & 1.2 & 1.2.1}$, and uses colors
to indicate the most recent move and the winning lines.

Modern students are able to see immediately that the tree above is
something that could be drawn by a computer interface, and it is easy
for them to imagine that the interface has buttons that turn colors on
and off, that turns arrows on and off, that switches between that
layout and Hutton's, and that displays or hides other parts of the
tree. {\sl The full tree would be so big that we don't want to see all
  of it} -- lots of objects in Mathematics are like that, and I will
call them {\sl expandable objects}.

\msk


{\bf Set comprehensions.} Set comprehensions are also expandable
objects. We can calculate the results of simple set comprehensions,
like $\setofst{10a}{a∈\{2,3,4\}}$, by ``reading them aloud'' in the
right way, i.e., by translating the mathematical notation into
English; but if we need to explain our mental process, or to handle
more complex cases, we need other techniques -- like 1) annotating the
parts of the comprehension as ``generators'', ``filters'', and
``resulting expression'', 2) converting the two standard notations for
comprehensions, $\setofst{10a}{a∈\{2,3,4\}}$ and
$\setofst{a∈\{2,3,4\}}{a^2<5}$, to a unified notation in which the
resulting expression always appears at the right, 3) converting them
to programs, and 4) using trees to calculate the result. The figure
below shows (1), (2), and (3):
%
\pu
%
\pu
%
\pu
$$\scalebox{0.8}{$
  \begin{array}{rcll}
  \ga{compr a^2<10} &=& \ga{compr a^2<10 u}  \\
                    &=& \ga{compr a^2<10 u;} & \ph{mm} \vbtsmash{compr a^2<10} \\
                    &=& \{4,9\}              \\
  \\
  \ga{compr 10*a}   &=& \ga{compr 10*a   u}  \\
                    &=& \ga{compr 10*a   u;} & \ph{mm} \vbtsmash{compr 10*a} \\
                    &=& \{20,30,40\}         \\
  \end{array}
  $}
$$

\unitlength=4pt
\def\closeddot{\circle*{0.6}}

Here is an example of (4):
%
%
$$\ga{compr triangle u;} \;=\;
  \ga{compr triangle dots} 
$$
$$\scalebox{0.8}{$
  \vbtsmash{compr triangle}
  \qquad
  \myvcenter{\ga{compr triangle tree thin}}
  \qquad
  \myvcenter{\ga{compr triangle tree wide}}
  $}
$$

Note that the two tree/tables above are ``simple'' in different ways:
in the one in the center it is obvious how we chose its columns, but
the one at the right is simpler to follow. See \cite{OchsNSC2026} for
more details and some class activities using set comprehensions.

\msk


{\bf Proofs with justifications.} For example:
%
$$\includegraphics[width=12cm]{2026logica-para-pessoas/example-6-just-4.pdf}$$
Imagine that this is a proof of $((6x^3)(7x^4))' = 294x^6$ that has a
computer interface. We could toggle off the display of the
intermediate steps, and we would get just $((6x^3)(7x^4))' = 294x^6$,
but we told the computer a) to display it as a series of steps that
are easy to justify, b) to number its equalities, c) to highlight what
changed between the left hand side and the right hand side of
`$\eqnp{4}$', d) to show the name of the rule that was used there e)
to show that rule ``as a formula'', f) to show what particular case we
used there ($n=3$), and g) to show what that formula becomes in that
particular case.


\section{Free and dependent variables}

Some people find dependent variables confusing. I was an extreme case
of that; I am not anymore, but I do believe that the students who find
dependent variables confusing deserve support -- like translations.


This is the first example of the chain rule in
\cite[p.101]{ApexCalculus4}, slightly reordered:
$$\def\Text#1{\text{#1:}\;\;\;}
  \begin{array}{lrcl}
  \Text{If}   & y     &=& f(g(x))        \\
  \Text{then} & y'    &=& f'(g(x))g'(x). \\
  \Text{If}   & f(x)  &=& x^2            \\
  \Text{and}  & g(x)  &=& 1-x            \\
  \Text{then} & f'(x) &=& 2x,            \\
              & g'(x) &=& -x,            \\
              & y     &=& (1-x)^2,       \\
              & y'    &=& 2(1-x)·(-1).   \\
  \end{array}
$$

We can use set comprehensions to make this more concrete to students
who don't understand equations well because they have problems with
free variables. The last two lines above say that the derivative of
the function with this graph
$$\setofst{(x,y)∈\R^2}{y=(1-x)^2}$$
is the function with this graph,

$$\setofst{(x,y')∈\R^2}{y'=2(1-x)·(-1)}$$
and we can change the $y'$ to $y$.

In the simplest examples, in which we can ignore the conditions on
smoothness and on domains, the chain rule is just
$\ddx f(g(x)) = f'(g(x))g'(x)$, and the derivative of the function
with this graph
$$\setofst{(x,y)∈\R^2}{y=f(g(x))}$$
is the function with this graph:
$$\begin{array}{cl}
    & \setofst{(x,y)∈\R^2}{y=\ddx f(g(x))} \\
  = & \setofst{(x,y)∈\R^2}{y=f'(g(x))g'(x)} \, .\\
  \end{array}
$$

Here is another way to write, and draw, that. If $f(x)=x^2$ and
$g(x)=1-x$ then $f'(x)=2x$ and $g'(x)=-1$, and:

$$\def\und#1#2{\underbrace{#1}_{\textstyle#2}}
  \def\draw#1{\includegraphics[width=2cm]{2026logica-para-pessoas/chainruleapex_00#1.pdf}}
  \sa {L1} {f(g(x))}
  \sa {L2} {\und{\ga{L1}}{(1-x)^2}}
  \sa {L3} {\und{\ga{L2}}{\draw{1}}}
  \sa {R1} {f'(g(x))g'(x)}
  \sa {R2} {\und{\ga{R1}}{2(1-x)·(-1)}}
  \sa {R3} {\und{\ga{R2}}{\draw{2}}}
  \scalebox{2.0}{$
    \ddx \ga{L3} \;=\; \ga{R3}
  $}
$$


\section{Maxima}

Maxima is a computer algebra system that is very extensible, and in
which we can define new operators and set lots of properties for any
operators -- for example we can set how they are \LaTeX{}ed, as in the
examples in the first column below. Note that \qco{matrix} is
standard, \qco{bmatrix} uses \qco{\\begin\{bmatrix\}} and
\qco{\\end\{bmatrix\}}, that use square brackets, \qco{barematrix}
doesn't draw any outer brackets, and \qco{verbatimmatrix} uses a kind
of verbatim box.

The internal representation of \co{f(a,b)} is \co{((\$F SIMP) \$A
  \$B)}, that is a tree in Lisp, in a format that is hard to
visualize. The result of \co{lisptreem(expr)} is the internal
representation of the result of \co{expr}, drawn as a 2D tree, as a
\qco{verbatimmatrix}.

The result of \co{lisptree2(expr)} includes the result of \co{expr}
and \co{lisptreem(expr)}.

The last input/output pair below, in \co{(\%i11)} and \co{(\%o11)},
shows how we can draw a series of equalities using a matrix.

\pu

\pu

\scalebox{0.6}{\def\colwidth{9cm}\firstcol{
  \maximagavbox{0cm}{9cm}{mintro matrixes 1}
}\def\colwidth{9cm}\anothercol{
  \maximagavbox{0cm}{9cm}{mintro matrixes 2}
}}

\newpage


Remember that in our Section 1 `ST'eve knew very few meanings for the
equality sign. We can use the code below to show him more meanings --
that are written in the same way in mathematical notation but
differently in Maxima:

\pu

\pu

\pu

\scalebox{0.8}{\def\colwidth{9cm}\firstcol{
  \maximagavbox{0cm}{9cm}{de-a-equacao-da-reta}
}\def\colwidth{9cm}\anothercol{
  \maximagavbox{0cm}{9cm}{de-a-equacao-da-reta 2}
}\def\colwidth{9cm}\anothercol{
}}

Typing \qco{2+3;} and then \co{RET} means ``calculate this and show
the result''; \qco{\$} means ``calculate this and don't show the
result''; \qco{:} is assignment; \qco{:=} defines a function; \qco{eq1
  : f(P[1]) = P[2]} stores an equality, as an expression, in \co{eq1};
\qco{solve([eq1,eq2], [a,b])} treats the equalities in \co{eq1} as two
equations to be solved and returns the solution as two other
equalities.

\def\und#1#2{\underbrace{#1}_{\textstyle#2}}

\msk


\subsection{Substitution}

This, that is an expanded version of the \co{(\%i7)} above,
$$\maximagavbox{0cm}{9cm}{de-a-equacao-da-reta 3}
$$
%
is similar to the operation that we were using to obtain particular
cases of equations, but here we are ``obtaing a particular case'' of
an expression, \qco{a*x+b}, that is not an equation. In the notation
that we used in section \ref{expandable-objects} this would be:
$$\und{(ax+b) \bmat{a:=-3 \\ b:=5}}{-3x+5}
$$
or:
$$(ax+b) \bmat{a:=-3 \\ b:=5} = -3x+5
$$

\msk


{\bf Juxtaposition.} In some cases, like $2(3+4)$, the juxtaposion
means an operation that was elided, but whose name we know, and we can
write it explicitly: $2·(3+4)$. In other cases the operation doesn't
have a standard notation, and we have to improvise. Let's use the
symbol `\standout{ap}' for application, and `\standout{s}' for
substitute:
$$\begin{array}{rcl}
  2(y+z) & \rarr & 2·(y+z) \\
  f(y+z) & \rarr & f \; \standout{ap} \; (y+z) \\
  (a+b)[a:=42] & \rarr & (a+b) \; \standout{s} \; [a:=42] \\
  \end{array}
$$

%
%
%
The operation `\standout{s}' is common in Logic and Computer Science
-- see for example \cite[p.6]{Harper} and
\cite[p.7]{HindleySeldin2008} -- but it is mostly treated as
``obvious'' in other places; some exceptions are \cite[p.146, p.149,
p.167]{Kindt} and \cite[p.482 and p.488]{FreudenthalDPh}.

\msk

\sa{s a by 42}{\text{$a$ by 42}}
\sa{s a by 42}{\mat{\text{substitute}\\
                    {\text{$a$ by 42}}}}


{\bf Process and object.} We can see the substitution as a process,
drawn as an arrow below,
$$a+b \;\; \xrightarrow{\ga{s a by 42}} \;\; 42+b
$$
%
%
that was reified, as in \cite[p.44]{Sfard}, and written as a suffix:
$$\sa{a by 42} {\text{$a$ by 42}}
  \sa{a:=42}   {\und{\mathstrut a:=42}{\ga{a by 42}}}
  \sa{[a:=42]} {\und{[ \; \ga{a:=42} \;]}
                    {\ga{s a by 42}}}
  (a+b) \;\; \ga{[a:=42]} = 42+b
$$

\msk


\pu

\pu

\pu

{\bf Infix.} In Maxima we can define new operators, and we can make
them infix, prefix, postfix, or n-ary. So we can do this, and
implement `\standout{s}' as \qco{\_s\_}:
$$\scalebox{1.0}{
    \maximagavbox{0cm}{9cm}{infix 1}
  }
$$


{\bf Simplification and lazyness.} The result in \qco{(\%o3)} above
was $b+42$, not $42+b$, because Maxima (almost) always simplifies
expressions, and its algorithm for simplifying sums can reorder and
collapse terms, and it can transform, for example, \co{2+3+4x+5x} into
\co{9x+5}.

Some functions in Maxima have versions that are less ``active'', in
the sense that they perform fewer simplifications. For example:
$$\scalebox{0.8}{$
    \maximagavbox{0cm}{9cm}{infix 2}
  $}
$$
In the terminology of Maxima \co{diff} is a ``verb'' and \co{'diff} is
a ``noun''.

\ssk

Another way to obtain operations that are simplified less, or that are
not simplified at all, is to create new operations and copy most of
the properties from the original versions to the copies, but omit the
properties that say how they are simplified. In the example below
\qco{+.} is a variant of \qco{+} that was created in that way:
$$\scalebox{0.8}{$
    \maximagavbox{0cm}{9cm}{infix 3}
  $}
$$

I call these less active operations ``lazy operations''.

\msk

\newpage


\pu

{\bf Prettifying.} The mathematical notation for parallel substitution
uses `:='s instead of `='s, puts each `:=' on a different line, and
uses square brackets. We can implement that with an operation \qco{V}
that pretty-prints the kinds of substitutions that \co{subst} expects,
like \co{[a=2,b=3]}, as \qco{bmatrix}es. The example below also shows
\qco{\_ss\_}, \qco{\_sss\_}, and \qco{\_ssu\_}, that use \qco{V} to
print substitutions in nice ways:

$$\scalebox{0.8}{$
    \maximagavbox{0cm}{9cm}{pretty subst}
  $}
$$

\msk

\newpage


\pu

{\bf Substituting functions.} Substitutions like
$\bsm{f(x):=x^2 \\ f'(x):=2x}$ are used a lot in Calculus, but they
are tricky to implement -- let's see why. We want to have this:
$$f(g(t)) \; \bmat{f(x):=g(g(x)) \\ g(x):=f(f(x))} \;=\; g(g(f(f(t))))$$
but look at this code:
$$\scalebox{0.8}{$
  \maximagavbox{0cm}{9cm}{tricky subst}
  $}
$$
\noindent both \qco{subst} and \qco{psubst}, that are Maxima
built-ins, don't do what we expect when we use \qco{S1}, but they work
well with \qco{S2}, in which we translated the functions to lambdas.

A good solution is to make \qco{\_s\_} translate functions to lambdas
when needed. Note that the lines \co{(\%i7)} and \co{(\%i8)} above
both yield $g(g(f(f(x))))$ -- that's because \qco{\_sss\_} calls this
more complex \qco{\_s\_}, that calls a function called
\qco{\_\_s\_fstolambdas} to perform this translation.


\sa {ggff small} {f\left(g\left(t\right)\right)\,\begin{bmatrix}f\left(x\right) := g\left(g\left(x\right)\right)\cr g\left(x\right) := f\left(f\left(x\right)\right)\cr \end{bmatrix}=g\left(g\left(f\left(f\left(t\right)\right)\right)\right)}
\sa {ggff big} {\underbrace{\underbrace{f\left(g\left(t\right)\right)}_{\textstyle \myvcenter{\vbtbgbox{\vbox{\vbthbox{f}\vbthbox{|}\vbthbox{g}\vbthbox{|}\vbthbox{t}}}}}}_{\textstyle \myvcenter{\vbtbgbox{\vbox{\vbthbox{ap\char95 \char95 .\ \ \ \ }\vbthbox{|\ \ \ |\ \ \ \ }\vbthbox{f\ \ \ ap\char95 \char95 .}\vbthbox{\ \ \ \ |\ \ \ |}\vbthbox{\ \ \ \ g\ \ \ t}}}}}\,\underbrace{\underbrace{\underbrace{\begin{bmatrix}f\left(x\right) := g\left(g\left(x\right)\right)\cr g\left(x\right) := f\left(f\left(x\right)\right)\cr \end{bmatrix}}_{\textstyle \begin{bmatrix}f := \lambda x.\,g\left(g\left(x\right)\right)\cr g := \lambda x.\,f\left(f\left(x\right)\right)\cr \end{bmatrix}}}_{\textstyle \begin{bmatrix}f := \lambda x.\,g\left(g\left(x\right)\right)\cr g := \lambda y.\,f\left(f\left(y\right)\right)\cr \end{bmatrix}}}_{\textstyle \begin{bmatrix}f := \myvcenter{\vbtbgbox{\vbox{\vbthbox{λ\char95 \char95 .}\vbthbox{|\ \ |}\vbthbox{x\ \ g}\vbthbox{\ \ \ |}\vbthbox{\ \ \ g}\vbthbox{\ \ \ |}\vbthbox{\ \ \ x}}}}\cr g := \myvcenter{\vbtbgbox{\vbox{\vbthbox{λ\char95 \char95 .}\vbthbox{|\ \ |}\vbthbox{y\ \ f}\vbthbox{\ \ \ |}\vbthbox{\ \ \ f}\vbthbox{\ \ \ |}\vbthbox{\ \ \ y}}}}\cr \end{bmatrix}}=\underbrace{\underbrace{g\left(g\left(f\left(f\left(t\right)\right)\right)\right)}_{\textstyle \myvcenter{\vbtbgbox{\vbox{\vbthbox{g}\vbthbox{|}\vbthbox{g}\vbthbox{|}\vbthbox{f}\vbthbox{|}\vbthbox{f}\vbthbox{|}\vbthbox{t}}}}}}_{\textstyle \myvcenter{\vbtbgbox{\vbox{\vbthbox{ap\char95 \char95 \char95 \char95 .\ \ \ \ \ \ }\vbthbox{|\ \ \ \ \ |\ \ \ \ \ \ }\vbthbox{λ\char95 \char95 .\ \ ap\char95 \char95 \char95 \char95 .}\vbthbox{|\ \ |\ \ |\ \ \ \ \ |}\vbthbox{x\ \ g\ \ λ\char95 \char95 .\ \ t}\vbthbox{\ \ \ |\ \ |\ \ |\ \ \ }\vbthbox{\ \ \ g\ \ y\ \ f\ \ \ }\vbthbox{\ \ \ |\ \ \ \ \ |\ \ \ }\vbthbox{\ \ \ x\ \ \ \ \ f\ \ \ }\vbthbox{\ \ \ \ \ \ \ \ \ |\ \ \ }\vbthbox{\ \ \ \ \ \ \ \ \ y\ \ \ }}}}}}

\pu

Here's a way to visualize how this substitution works:
$$\ga{ggff big}$$

The path is $\psm{•&&•\\↓&&↑\\•&→&•}$, and in the last part, the
`$↑$', Maxima handles the $β$-reductions in the standard way, in which
Church-Rosser theorem guarantees that the order of the reductions
doesn't matter. See \cite[pages 7, 11, and 14]{HindleySeldin2008} for
the details, and note that the clause (g) in p.7 ``{\sl chang[es]
  bound variables to avoid clashes}''.

\msk

A few years ago I thought that the right thing to do would be to
implement this `[:=]' in Lean ``to show how natural it is''. That
turned to be much harder than I expected, and to make things worse
Lean doesn't have a REPL and doesn't have a built-in \LaTeX{}
generator... while in Maxima \qco{\_\_s\_fstolambdas} is just 9 lines
and the rest of \qco{\_s\_} is just 6 lines. I wanted to convince
people that this `[:=]' is a natural operation, {\sl even though I
  couldn't find it explained in details in the literature}, and
showing how it can be implemented in 15 lines in Maxima seemed to be a
good way.

I will explain more reasons for not using Lean in the next section.


\newpage


\section{Why not Lean?}

The ``right way'' to handle ``procedural students'' that don't know
anything about Logic would be to teach Logic to them. There is a very
good example of how to do that in \cite{Yalep} -- note in particular
the ``Pedagogical progression'' in its Appendix B -- and the authors
have a comparison with other computational approaches in
\cite{YalepSurvey}... {\sl but we don't have time for the right way.}

Let me quote a few sentences from
\cite[p.30]{BottoniTFM}:

\begin{quote}

  However, Lean is not currently suitable for high school or
  undergraduate university teaching. High school mathematics curricula
  in Switzerland are already crowded, and introducing abstract topics
  like first-order logic or natural induction would push the limits.

\end{quote}

The syllabus of Calculus 2 in my campus is huge, and it was hard to
cover it completely even when the students came in much better
prepared. I needed an approach in which a) the non-textbook material
would be just a small part of each class, b) the students with little
experience with computers would be able to follow almost everything,
and c) the programs would be easy to install at home, and easy to use
both at home and at the computer lab. So I ended up using Maxima with
the interface described in \cite{OchsEmacsConf2024}; my presentation
at the EBL 2025 (\cite{OchsEBL2025}) was about how people who
installed Maxima in that way could install Lean easily and then follow
a basic tutorial on Lean, but in Calculus 2 I show Lean very briefly,
just as a curiosity.

\pu

\pu
\section{Holes}

Maxima comes with a function called \qco{dpath}, that expects an
object and a ``path'', and `d'isplays a box around subobject obtained
by following that path. In the examples below the path is \qco{2,1},
that means ``second child node, then first child node'':
$$\scalebox{0.8}{$
  \maximagavbox{0cm}{9cm}{holes 1}
  $}
$$

In the last line we used \qco{substpart("?", obj, <path>)}, that does
something similar to \qco{dpart} but replaces the object there by a
\qco{"?"} instead of by a box.


Note that in the example above we started with $2+3+4=2+7$, that is
{\sl true}, and we obtained $2+3+4=\Rq+7$, that is a {\sl problem}, or
an {\sl exercise}. I will refer to the `$\Rq$' as a {\sl hole}; see
\cite[p.465]{FreudenthalDPh}.


Many exercises in Calculus books can be rewritten as true expressions
with holes. For example, the exercise 11 in
\cite[p.109]{ApexCalculus4} is ``compute the derivative
$\ddx (\ln x + x^2)^3$ using the chain rule'', and we can rewrite it
as a true expression with holes by doing this:
$$\scalebox{0.8}{$
  \maximagavbox{0cm}{14cm}{holes 2}
  $}
$$

This program is a {\sl manipulable object}: students can, for example,
duplicate the 4-line command with the definition of \qco{S}, modify
the substitutions in the copy, and run the new code in the REPL to see
what changes in the result...

We can also ask students to modify the \qco{S} above to obtain true expressions
that solve other exercises of the book.

\msk

Very few students understand right away how they can find the
derivative of $\ddx \sin(42x)$ by this method by solving a simpler
problem first -- the simpler problem is
$\ddx f(g(x)) \bsm {f(x):=\Rq \\ g(x):=\Rq} = \ddx \sin(42x)$. By
making the students work together the weaker students will see that
the stronger students found a way to decompose bigger problems into
smaller problems, and will see that there is a method there that they
can learn.

Note that I haven't shown (yet!) how to transform our expandable proof
from section \ref{introduction} into a manipulable object -- I am just
showing how to manipulate objects that look like its {\sl
  justifications}.


\section{Extremely Bad Foundations}
\label{bad-foundations}


In \cite{FeynmanJoking} Richard Feynman tells several stories about
his stays in Brazil and his experiences teaching classes here. This is
from one of them, about a class that he taught in 1952:

\begin{quote}

  I taught a course at the engineering school on mathematical methods
  in physics, in which I tried to show how to solve problems by trial
  and error. It's something that people don't usually learn, so I
  began with some simple examples of arithmetic to illustrate the
  method. I was surprised that only about eight out of the eighty or
  so students turned in the first assignment. So I gave a strong
  lecture about having to actually {\sl try} it, not just sit back and
  watch {\sl me} do it.

  After the lecture some students came up to me in a little
  delegation, and told me that I didn't understand the backgrounds
  that they have, that they can study without doing the problems, that
  they have already learned arithmetic, and that this stuff was
  beneath them.

\end{quote}


We use the term ``Foundations of Mathematics'' for the axioms and
logical rules behind the maths that a person does; nowadays the most
common foundations are ZFC and Type Theories. I will use the
(non-standard) term ``Bad Foundations'' for people who don't know any
logical rules and who do mathematics by just memorizing methods and
applying them, and ``Extremely Bad Foundations'' for people who not
only have bad foundations but who also believe that all methods that
require trial and error -- or even trial-and-improve; see
\cite[p.180]{DrijversBoonReeuwijk} -- are wrong, and they refuse to
learn them.

\msk

Some courses are supposed to work like sieves that only let through
the students whose mathematical skills are fine enough. But the system
was already dysfunctional, and then we had COVID and now LLMs, and now
lots of very coarse students with Extremely Bad Foundations are
reaching classes like Calculus 2, and we need to do some something for
them that is better than just failing them, or, even worse, letting
them pass these advanced courses unchanged...

\msk


One idea -- {\sl that I am treating as a work in progress; I have only
  a few examples ready, and they were tested only on small classes,
  and without a rigorous methodology} -- is that we can use Computer
Algebra Systems like Maxima to treat some methods from Calculus as
expandable objects, make the computer draw the expanded versions of
some examples, and use that to show to students with extremely bad
foundation how the other students think, and what are the steps that
they don't write down. It should be possible to use ``holes'' to
convert these expanded version into exercises, and if the students
with bad and extremely bad foundations do these exercises then it {\sl
  may be possible} that they will develop certain pattern-matching
skills by drill, and they {\sl may} reach the point that John Mason
describes in \cite[p.4]{MasonShift}, in which there happens ``{\sl a
  splitting of attention, [with] one part [of the mind] involved in
  calculation while another part remains aloof but observant}''. After
learning certain new procedural skills they {\sl may} have an ``aha!
moment'' and start to understand how Logic works.


\msk

I will close with another example of a method as a manipulable object.
\cite[p.639]{Stewart8} explains ODEs with separable variables starting
by the general case, $\dydx=\frac{g(x)}{h(y)}$, then solving it by a
method that involves differentials and includes steps that many
students find hard to believe; and then its Example 1, in the
following page, is $\dydx = \frac{x^2}{y^2}$, that has solutions that
are very hard to draw by hand...


\sa {EDOVS def [M]} {\ga{[M]}=\begin{pmatrix}{\frac{\textit{dy}}{\textit{dx}}}&\mbox{ = }&{\frac{g\left(x\right)}{h\left(y\right)}}\cr h\left(y\right)\,\textit{dy}&\mbox{ = }&g\left(x\right)\,\textit{dx}\cr \int {h\left(y\right)}{\;dy}\big.&\mbox{ = }&\int {g\left(x\right)}{\;dx}\big.\cr \mbox{ || }&&\mbox{ || }\cr H\left(y\right)+C_1 &\mbox{ = }&G\left(x\right)+C_2 \cr H\left(y\right)&\mbox{ = }&G\left(x\right)+C_3 \cr H^{-1}\left(H\left(y\right)\right)&\mbox{ = }&H^{-1}\left(G\left(x\right)+C_3 \right)\cr \mbox{ || }&&\cr y&&\cr \end{pmatrix}}
\sa {EDOVS def [A]} {\ga{[A]}=\begin{pmatrix}{\frac{\textit{dy}}{\textit{dx}}}&\mbox{ = }&{\frac{-\left(2\,x\right)}{2\,y}}\cr 2\,y\,\textit{dy}&\mbox{ = }&-\left(2\,x\right)\,\textit{dx}\cr \int {2\,y}{\;dy}\big.&\mbox{ = }&\int {-\left(2\,x\right)}{\;dx}\big.\cr \mbox{ || }&&\mbox{ || }\cr y^2+C_1 &\mbox{ = }&-x^2+C_2 \cr y^2&\mbox{ = }&-x^2+C_3 \cr \sqrt{y^2}&\mbox{ = }&\sqrt{-x^2+C_3 }\cr \mbox{ || }&&\cr y&&\cr \end{pmatrix}}
\sa {EDOVS def [S5]} {\ga{[S5]}=\begin{bmatrix}g\left(x\right) := -\left(2\,x\right)\cr h\left(y\right) := 2\,y\cr G\left(x\right) := -x^2\cr H\left(y\right) := y^2\cr H^{-1}\left(y\right) := \sqrt{y}\cr \end{bmatrix}}

\sa{[M]}{\CFname{M}{}}
\sa{[A]}{\CFname{A}{}}
\sa{[S5]}{\CFname{S5}{}}

For me the ``archetypal'' ODE with separable variables is
$\dydx = -\frac{x}{y}$, whose solutions are circles centered on
$(0,0)$, and the general method can be reconstruced from that
archetypal case. The archetypal case and the general method are:
$$\scalebox{0.8}{$
  \def\mbox#1{#1}   
  \ga{EDOVS def [A]}
  \quad
  \ga{EDOVS def [M]}
  $}
$$
and the techniques for working on the ``archetypal case'' and the
``general case'' in parallel are explained in \cite{OchsIDARCT} and
\cite{OchsMD}. {\sl I hope to finish the details for this example
  soon.}





\printbibliography

\pu

\GenericWarning{Success:}{Success!!!}  

\end{document}